\documentclass[graybox]{svmult}

\usepackage{amsmath}
\usepackage{amsthm}
\usepackage{mathtools}

\usepackage{type1cm}        
\usepackage{makeidx}         
\usepackage{graphicx}        
\usepackage{multicol}        
\usepackage[bottom]{footmisc}

\usepackage{subfigure}

\usepackage{newtxtext}       %
\usepackage[varvw]{newtxmath}       

\usepackage{bbm}
\usepackage[normalem]{ulem}
\usepackage{hyperref}
\usepackage{cleveref}

\newtheorem{algorithm}{Algorithm}
\newtheorem{claimp}{Claim}

\makeindex             

\begin{document}

\title*{Studying wildfire fronts using advection--diffusion--reaction models}
\titlerunning{Studying wildfire fronts using advection--diffusion--reaction models}
\author{Koondanibha Mitra\orcidID{0000-0002-8264-5982},
Qiyao Peng\orcidID{0000-0002-7077-0727} and\\
Cordula Reisch\orcidID{0000-0003-1442-1474}
}
\institute{Koondanibha Mitra \at Mathematics \& Computer Science Department, Eindhoven University of Technology, 5600 MB Eindhoven, The Netherlands, \email{k.mitra@tue.nl}
\and Qiyao Peng \at Mathematical Institute, Faculty of Science, Leiden University, Neils Borhweg 1, 2333 CA, Leiden, The Netherlands. \email{q.peng@math.leidenuniv.nl} 
\and Cordula Reisch \at Institute for Partial Differential Equations, TU Braunschweig, Universitaetsplatz 2, 38106 Braunschweig, Germany  \email{c.reisch@tu-braunschweig.de}}
%
%
\maketitle

\abstract*{Each chapter should be preceded by an abstract (no more than 200 words) that summarizes the content. The abstract will appear \textit{online} at \url{www.SpringerLink.com} and be available with unrestricted access. This allows unregistered users to read the abstract as a teaser for the complete chapter.
Please use the 'starred' version of the \texttt{abstract} command for typesetting the text of the online abstracts (cf. source file of this chapter template \texttt{abstract}) and include them with the source files of your manuscript. Use the plain \texttt{abstract} command if the abstract is also to appear in the printed version of the book.}
\vspace{-5em}
\abstract{In this work, we study the propagation of wildfires using an advection--diffusion--reaction model which also includes convective and radiative heat loss. An existing model is discussed \cite{asensio_2002} and a physically consistent modification of the model is proposed. Using this, the existence of travelling waves (TWs) in the one-dimensional case is investigated. Prior numerical studies reveal the existence of TWs \cite{reisch_2023}. Under the travelling wave ansatz and certain approximation, the model is reduced to a semi-autonomous dynamical system with three unknowns which can be analyzed by a shooting algorithm. It is hypothesized that under mild wind speeds, TWs in both directions exist, and under strong tailwinds only TWs in the direction of wind are possible. The theoretical implications are investigated using both solvers for the PDE models and the shooting algorithm. The results match, and unveil the dependence of the fronts on the parameters consistent with the predictions.}

\vspace{-2em}
\section{Introduction}\label{sec:intro}
\vspace{-1em}
Wildfire as a natural phenomenon poses a palpable threat to our society which has become evident in recent years. Due to climatic change, the frequency and severity of wildfires will only increase over the next decades \cite{senande_2022}. This makes it important to develop accurate mathematical models and analytical tools that predict how wildfires spread.
 For instance, wildfires are known to possess a front that moves with a given speed which depends heavily on the wind speed and flammability of the vegetation. Naturally, the question arises if its spread can be approximated by a \emph{travelling wave}, abbreviated as TW henceforth. To investigate this question,
 in this work, we focus on continuous advection--diffusion--reaction (ADR) models for wildfires which were studied for example in \cite{ asensio_2002, reisch_2023, burger_2020b}.
The ADR models describe the fire dynamics by modeling the variation of temperature and biomass. 
There exists a two-way and nonlinear coupling between these two primary variables. In addition, the reaction term is discontinuous only playing a role in the domain that is burning. This leads to challenges in analyzing the model and predicting the existence of travelling waves. Having estimates about TW solutions in terms of parameters such as wind speed will also help us to predict the areas at risk for planning interventions and potentially necessary evacuations.


Advection--diffusion--reaction (ADR) models for wildfires have been intensively studied: In \cite{weber1991} TW solutions for a diffusive model for the temperature in solid and gaseous phases are analyzed.  \cite{VolpertTW1994} gives an overview of TWs for static problems and for problems where both variables diffuse. 
\cite{weber1997} analyses the travelling wave speed for a fully diffusive combustion model with a continuous reaction function including the Arrhenius law. A major observation of the analytical and numerical investigations is that the spreading speed of the fire depends mostly on the fire switching `on' and `off' and not on the cooling down of the temperature. 
\cite{asensio_2002} proves the existence and uniqueness of weak solutions for the model including an unsteady switching function for the combustion process. Based on the existence results, a mixed finite element method is proposed. 
In \cite{mandel2008} the discontinuous combustion function is approximated by a continuous function including the Arrhenius law and linear cooling processes. The approximation function is chosen such that TW solutions can be observed numerically. Data assimilation is discussed for estimating the model parameters. 
\cite{burger_2020b, burger_2020} propose implicit-explicit numerical methods for simulating the ADR model.  The spread of fire is numerically studied in heterogeneous environments and under the influence of wind. 

Besides the numerical studies and the analytical proof of the existence and uniqueness of weak solutions, results also indicate the existence of TWs for models with one non-diffusive variable.
In \cite{babak_2009}, the existence of TW solutions was proven for a simplified model where ambient cooling effects were not regarded. In this simplified model, there is no mechanism by which the fire can die out and hence the temperature does not decrease. This leads to an unphysical situation. 
Instead, in the current work, inspired by \cite{MitraTW2023} where TWs having similar profiles were analyzed for the spreading of biofilms, we make several important steps toward showing the existence of TWs in a much more general setting. In \Cref{sec:ADRmodel}, we discuss an existing ADR model and propose its modification. Then we approximate the nonlinear combustion term and non-dimensionalize the system. In \Cref{sec:TW}, a TW transformation is proposed for the one-dimensional case leading to an ordinary differential equation system in the wave variables. 
In \Cref{sec:Num}, a shooting algorithm is proposed to find the TW solutions directly for the approximate model. Additionally, the dynamics of the partial and the ordinary differential equation systems are numerically compared to verify the transformation and our predictions.


\section{Advection--diffusion--reaction models for wildfire}\label{sec:ADRmodel}
\vspace{-1em}
\subsection{The existing ADR model}
\vspace{-1em}
The model proposed in \cite{asensio_2002} and studied in \cite{reisch_2023} reads
 \begin{equation}  \label{eq:original_model}
     \left\{ \begin{aligned}
       \rho_0 C  \left(\frac{\partial T}{\partial t}  + \mathbf{w}\cdot\nabla T \right) & =  \nabla\cdot \left[ K(T) \nabla T \right]     - h(T-T_{\infty}) +  \Psi(t)\rho_{0} H Y,   \\
 \frac{\partial Y}{\partial  t} &= -\Psi(t)Y,
     \end{aligned}\right.
 \end{equation}
with the temperature $T$ and the biomass $Y>0$. The model parameters are the bulk density $\rho_0$, the specific heat $C$, the convection coefficient $h$, the ambient air temperature $T_\infty$, the heating value $H$, and the wind speed $\mathbf{w}$. The first term on the right denotes the combined effect of heat diffusion and radiation, with the resulting diffusivity being
\begin{equation*}
    K(T)= k + \epsilon\sigma T^3, k:\text{heat diffusivity}, \,\sigma:\text{Stefan-Boltzman constant,}\, \epsilon:\text{scaling factor.} 
\end{equation*}
See \cite{asensio_2002} for the details.
The second term is convective heat loss with the convection constant $h>0$, and the last term denotes the heat generated by combustion. The function $\Psi$ describes this nonlinear process (with $\mathbbm{1}(B)$ denoting characteristic function of a measurable set $B$)
 \begin{equation}\label{eq:Psi_1st}
 \Psi(t) =\Psi_0(T(t)):=A \mathbbm{1}(\{T(t)\geqslant \Bar{T}\})  \exp\left(-\frac{T_{\rm ac}}{T}\right)
 \end{equation}
 with the activation temperature $T_{\rm ac}$, the ignition temperature $\Bar{T}$ and a pre-exponential factor $A>0$ arising from the Arrhenius formula. Model \eqref{eq:Psi_1st} will henceforth be called the \emph{extinguishing fire} model or \emph{EF} since after ignition, the fire stops if $T<\Bar{T}$. Suitable boundary conditions and initial conditions (with maximum biomass $Y_{\rm ref}>0$) complete the problem \eqref{eq:original_model}--\eqref{eq:Psi_1st}.

 \vspace{-1.5em}

 \subsection{Model variations}
\vspace{-1em}

In this work, we propose a slight variation of the expression \eqref{eq:Psi_1st} for $\Psi$ which we believe to be more realistic, and which helps in the travelling wave analysis later on:
\begin{align}\label{eq:new_psi2}
    \Psi(t)=\Psi_2(\{T(\tau): \tau\in [0,t]\}):=A \mathbbm{1}\left(\bigcup_{0\leqslant\tau\leqslant t}\{T(\tau)\geqslant \bar{T} \}\right) \exp\left(-\frac{T_{\rm ac}}{T}\right).
\end{align}
This new expression assumes that the fire remains turned `on' once it has started at a location (after the temperature exceeds $\bar{T}$) which is in contrast to \eqref{eq:Psi_1st} which assumes that it immediately gets extinguished if $T<\Bar{T}$. Hence, it will be called the \emph{persistent fire} model or \emph{PF} moving forward. This is more physical in our opinion and is supported by aerial observations of wildfires.
Furthermore, \cite{weber1997} assesses the higher importance of the on-setting of combustion to the travelling wave speed compared to the cooling processes at the end of the burning (tail region). 

\begin{figure}
    \centering
    \includegraphics[width=.5\textwidth]{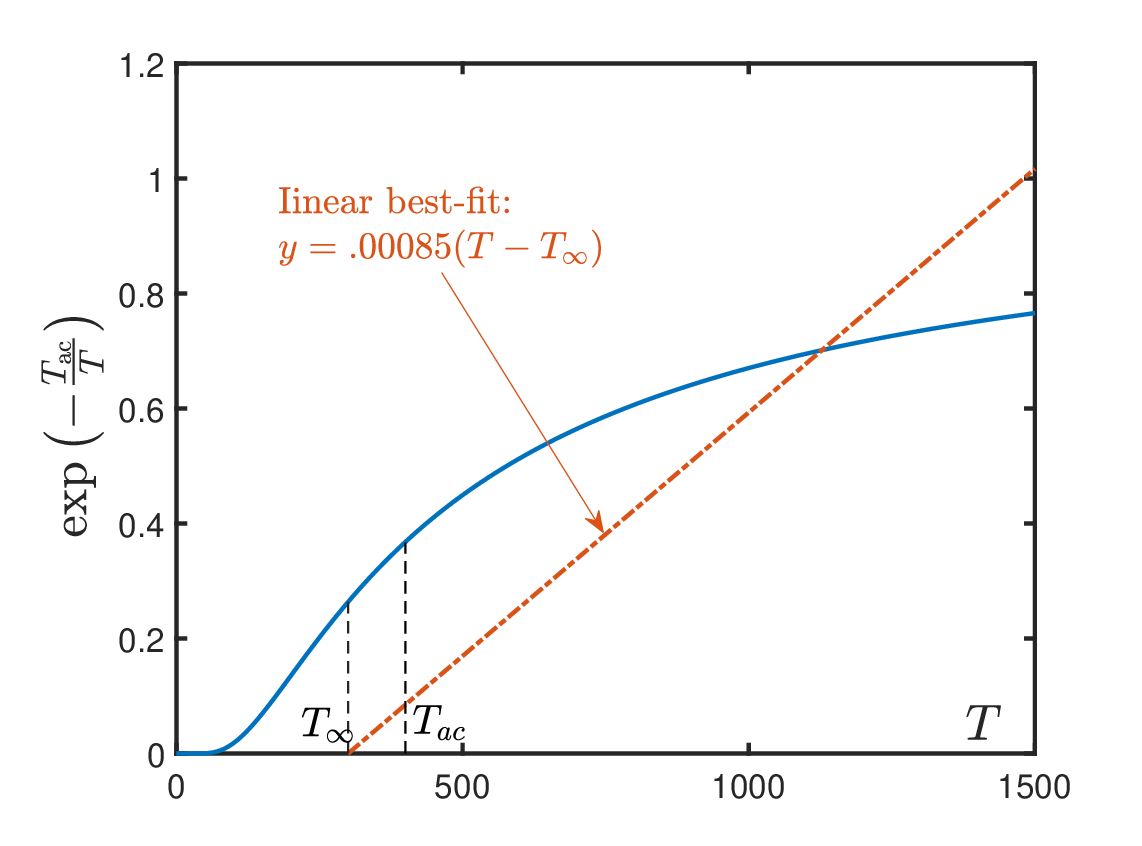}
    \caption{The exponential function $\exp(-T_{\rm ac}/T)$ and its least-square fitted linear approximation in the interval $300\leq T\leq 1500$. Here $T_{\rm ac}=400$.}
    \label{fig:linear_approx}
    \vspace{-2em}
\end{figure}

Purely for the sake of analysis, we would further assume
\begin{equation}\label{eq:linear_approx}
    \exp\left(-\frac{T_{\rm ac}}{T}\right)\approx \Lambda (T-T_\infty)
\end{equation}
for a constant $\Lambda>0$ that can be obtained by least-square fitting in the physically relevant temperature range, see \Cref{fig:linear_approx}. This linearizes the exponential term, making algebraic manipulation possible. Moreover, it introduces a relatively small error in the intermediate temperature range as $\exp(-T_{\rm ac}/T_\infty)$ is small which implies that the zeroth order term can be ignored.

  \vspace{-2em}

 \subsection{Non-dimensionalization and effective model}
\vspace{-1em}
Observing that the diffusion, convection, and heat production terms have to be roughly in balance for the formation and propagation of the wildfire, we propose the following non-dimensionalization. Define the scaled variables
\begin{align}
  &\tilde{t}:=\frac{t}{t_{\rm ref}},  \quad \tilde{\mathbf{x}}:=\frac{\mathbf{x}}{L_{\rm ref}},\quad u:= \frac{T-T_{\infty}}{\bar{T}-T_\infty},\quad v:=\frac{Y}{Y_{\rm ref}},\end{align}
  with reference time $ t_{\rm ref}:= \frac{C}{A\Lambda H}$ and length $L_{\rm ref}:= \sqrt{\tfrac{k}{\rho_0 A \Lambda H}}$ resulting in the dimensionless quantities
  \begin{align}
  \tilde{h}:=\frac{h}{\rho_0 A \Lambda  H},\quad \tilde{\mathbf{w}}:=\frac{t_{\rm ref}\,\mathbf{w}}{L_{\rm ref}},\quad \gamma:=\frac{C(\Bar{T}-T_\infty)}{H}, \quad \tilde{K}(u(T)):=\frac{K(T)}{k}.
\end{align}
With the abuse of notation, dropping the $\tilde{}$ in the quantities above, we get from \eqref{eq:original_model} with \eqref{eq:new_psi2}--\eqref{eq:linear_approx}, the ADR-PF model with linear approximation:
 \begin{align}  \label{eq:model_uv}
     \left\{ \begin{aligned}
     	\partial_t u + \mathbf{w}\cdot\nabla u &= \nabla \cdot[K(u) \nabla u] -h\,u + u\,v \,\mathbbm{1}(\cup_{0\leqslant\tau\leqslant t}\{u(\tau)\geqslant 1\}) \\
	\partial_t v &= - \gamma\, u\,v\,  \mathbbm{1}(\cup_{0\leqslant\tau\leqslant t}\{u(\tau)\geqslant 1\} )
          \end{aligned}\right.
\end{align}


\vspace{-2em}
\section{One-dimensional fronts and travelling waves}\label{sec:TW}
\vspace{-1em}
In this section, we restrict our discussion to the one-dimensional case. Hence, $\nabla$ in \eqref{eq:model_uv} is replaced by $\partial_x$. Without loss of generality, assume wind speed $w\geq 0$.  For analyzing TWs, the natural boundary condition is to impose that the temperature reaches the atmospheric temperature as $x\to \pm \infty$, and as a consequence, the diffusive/radiative heat fluxes vanish. On the other hand, ahead of the fire, the biomass remains constant at $Y_{\rm ref}$. In dimensionless terms, the boundary and initial conditions satisfy
\begin{subequations}\label{eq:ICBC}
   \begin{align}
    &u(\pm \infty,t)=\partial_x u(\pm \infty,t)=0 &&\forall\, t>0, \\
    &u(x,0)< 1, \quad v(x,0)= 1 \quad &&\forall\, x> 0
 \end{align}
with $u(0,0)=1$ used to fix the fire-front.
\end{subequations}



 \noindent
 The \textbf{travelling wave ansatz} in this case is the assumption that solutions of the form 
 \begin{align}
 u(x,t)= u(\xi(x,t)), \;\; v(x,t) = v(\xi(x,t))   
 \end{align}
 exists for the travelling wave coordinate $\xi=x-ct$ with an unknown speed $c\not=0$. 
The TW ansatz transforms System \eqref{eq:model_uv} into 
\begin{align}\label{eq:transformed_system}
\left\{		\begin{aligned}
    (w-c) \frac{\mathrm{d} u }{\mathrm{d} \xi } & = \frac{\mathrm{d}}{\mathrm{d} \xi} \left( K(u) \frac{\mathrm{d} u}{\mathrm{d} \xi} \right) -h\,u + u\,v\, \mathbbm{1}(\{\xi\leq  0\}) , \\
		c \frac{\mathrm{d} v}{\mathrm{d} \xi} & =  \gamma\, u\,v\,  \mathbbm{1}(\{\xi\leq 0\}).
\end{aligned}\right.
\end{align}
with the transformed boundary conditions consistent with \eqref{eq:ICBC}
\begin{align}\label{eq:bc_transformed}
	 u(\pm \infty) =\frac{\mathrm{d} u}{\mathrm{d} \xi} (\pm \infty) = 0, \quad u(\xi)<1, \quad v(\xi) =1\quad \forall\, \xi > 0, \quad\text{ and }\; u(0)=1. 
\end{align}
Since the above conditions imply that the fire has been on only in $\xi\leq 0$, the characteristic function in \eqref{eq:model_uv} simply reduces to $\mathbbm{1}(\{\xi\leq 0\})$ in this case.

To reduce \eqref{eq:transformed_system} further, the total thermal energy is introduced:
 \begin{align}
 	z = \frac{1}{c} \int^{\infty}_\xi u(\tilde{\xi}) \, \mathrm{d} \tilde{\xi}, \quad \Rightarrow  \quad c \frac{\mathrm{d} z}{\mathrm{d} \xi} =-u.
 \end{align}
Inserting $z$ and the expression of $\frac{\mathrm{d} v}{\mathrm{d} \xi}$ in \eqref{eq:transformed_system} and integrating in $\xi$, one finds the first integral (conserved quantity)
 \begin{align}\label{eq:first_integral}
 K(u) \frac{\mathrm{d} u}{\mathrm{d}\xi} + (c-w)u+ c\, h\, z + \frac{c}{\gamma} v = D \text{ (constant independent of $\xi$)}.
 \end{align}
 Passing the limit $\xi\to \pm\infty$, and using the boundary conditions \eqref{eq:bc_transformed} one obtains 
\begin{align}
    v(-\infty)=:r, \quad z(-\infty)=\frac{1}{\gamma h}(1-r) \quad \text{ for some } r\in (0,1).
\end{align}
The $\xi\to \infty$ limit further yields $D=\frac{c}{\gamma}$. Using the above, rearrangement of \eqref{eq:first_integral} yields the dynamical system for the triplet $(u,v,z)$ with unknown wave-speed $c\not=0$, 
\begin{equation}
    \label{eq:uvz_ODE_xi}
    \left\{
    \begin{aligned}
        K(u)\frac{du}{d\xi} &= (w-c)u-c\,h\,z+\frac{c}{\gamma}(1-v),\\
        c\frac{dv}{d\xi} &= \gamma\,u\,v\mathbbm{1}(\{\xi\leq 0\}),\\
        c\frac{dz}{d\xi} &= -u.
    \end{aligned}
    \right.
\end{equation}
The conditions \eqref{eq:bc_transformed} demand that for an unknown $r\in (0,1)$,
\begin{align}\label{eq:BC_ode}
    (u,v,z)\to \begin{cases}
        (0,1,0) &\text{ as } \;\xi\to +\infty,\\
        \left(0,r,\tfrac{1}{\gamma h}(1-r)\right) &\text{ as } \; \xi\to -\infty,
    \end{cases} 
\end{align}
with the $\xi$-cordinates being shifted so that the supremum of the set $\{\tilde{\xi}: u(\Tilde{\xi})>1\}$ is 0.
The points on the right of \eqref{eq:BC_ode} are the equilibrium points for the semi-autonomous system \eqref{eq:uvz_ODE_xi}. 

The above system can be investigated numerically as a shooting problem, i.e. solving \eqref{eq:uvz_ODE_xi} backward starting from $(0,1,0)$ for different values of $c$ to find when the condition at $\xi\to -\infty$ is satisfied. This is precisely what is done in the next section, see \Cref{algo:shooting}. Moreover, in a future work, we also conclude the existence of TWs by analyzing \eqref{eq:uvz_ODE_xi}--\eqref{eq:BC_ode}. Here, we summarize our claim without proof:

\begin{claimp}[Existence of travelling waves for the ADR-PF model] Let $h<1$ (necessary condition). Then under certain parametric constraints involving $h$ and $\gamma$, there exists $\bar{w}>0$ such that for all $w\in [0,\bar{w})$, at least two travelling wave solutions exist of the System \eqref{eq:uvz_ODE_xi}--\eqref{eq:BC_ode} with $c=c_+>0$ and $c=c_-<0$ respectively (with $c_++ c_->0$), whereas for $w>\bar{w}$, there exist only travelling waves with $c=c_+>0$.
\end{claimp}

The duplicity of the TWs results from the directional symmetry of the system. For example for $w=0$, the problem is completely symmetric with respect to the transformation $x\mapsto-x$, see also the numerical results in \cite{reisch_2023}. Moreover, we stress that no TWs exist if $h>1$, see (g)-(h) plots of \Cref{fig:PDE_Eq3_exp}.

\vspace{-2em}

\section{Numerical results}\label{sec:Num}
\vspace{-1EM}
In this study, the finite-element method (FEM) together with the Strang splitting is utilized for numerical simulation of the EF and the PF models. The dimensional values of the parameters used in the simulation are taken from Table 1 of \cite{reisch_2023}, i.e., $\rho_0=40kg/m^3$, $C=1kJ/(kg\cdot K)$, $k=2kW/(m\cdot K)$, $h=4kW/(m^3\cdot K)$, $A=0.05s^{-1}$, $H=4000kJ/kg$, and $T_{\rm ac}=400K$. The computational domain is $(-250, 250)m$ and the time step size is $\Delta t = 0.4s$. For the initial condition, the temperature in $(-25,25)m$ is $470K$, which is higher than the ignition temperature $\displaystyle\bar{T}=400K$, while the rest of the domain is in ambient temperature $T_\infty=300K$. In the beginning, the fuel concentration $Y$ is $1$ over the entire domain.

\begin{figure}[h]
    \centering
    \subfigure[No wind and TWs exist]{
    \includegraphics[width = 0.49\textwidth,trim={5 0 5 0},clip]{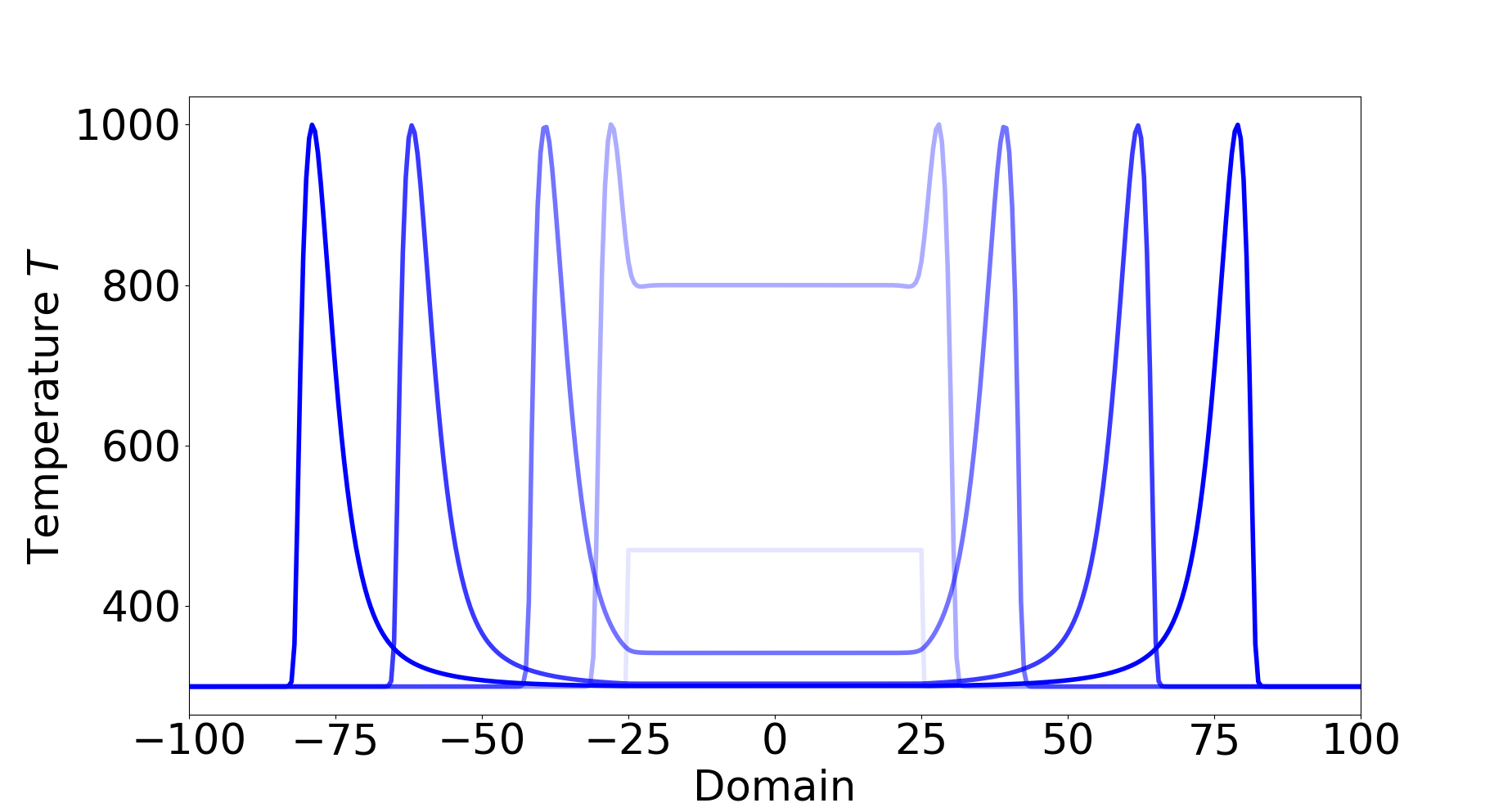}}
    \subfigure[No wind and TWs exist]{
    \includegraphics[width = 0.49\textwidth]{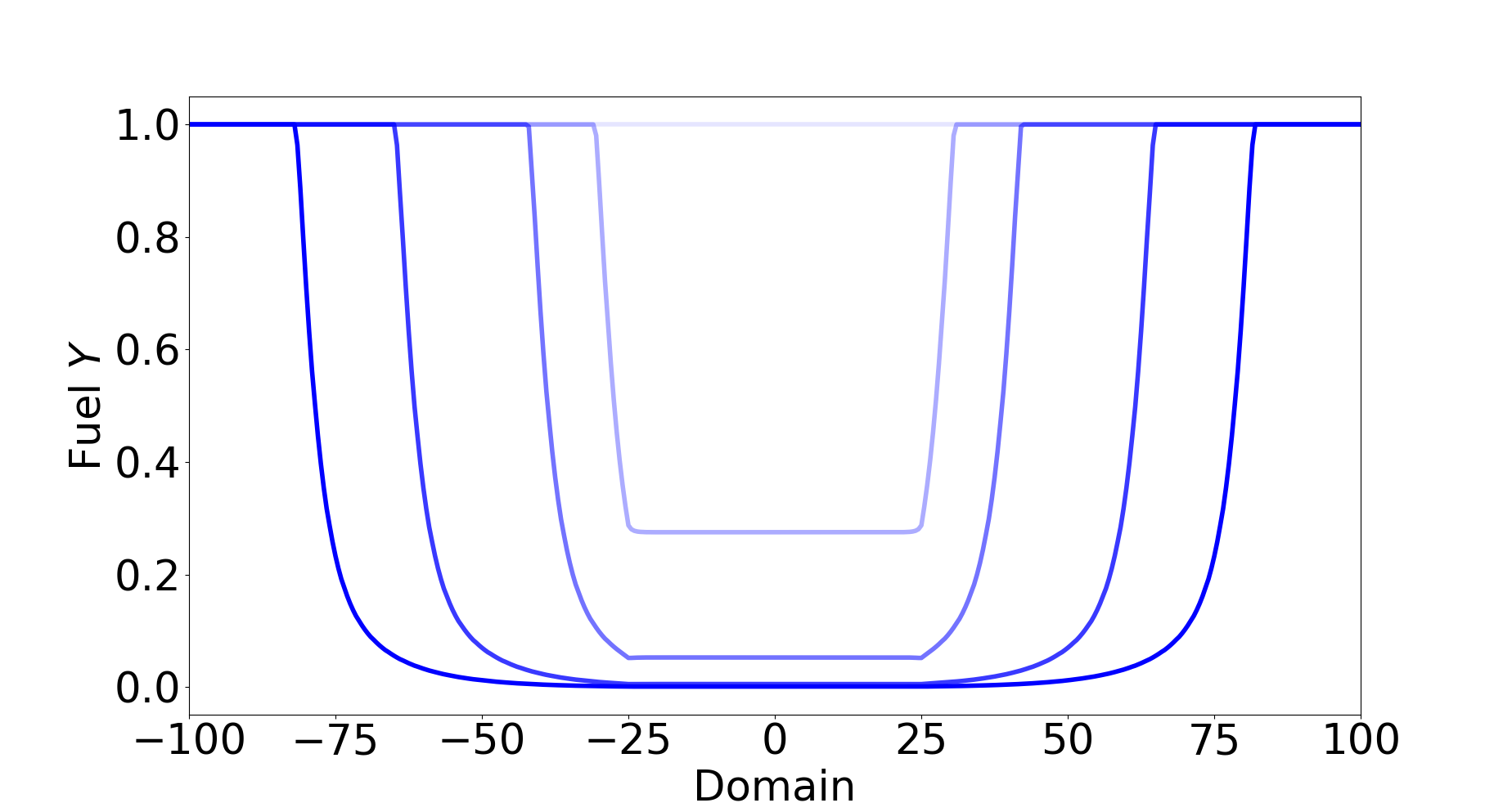}}
    \subfigure[Mild wind and TWs exist]{
    \includegraphics[width = 0.49\textwidth,trim={5 0 5 0},clip]{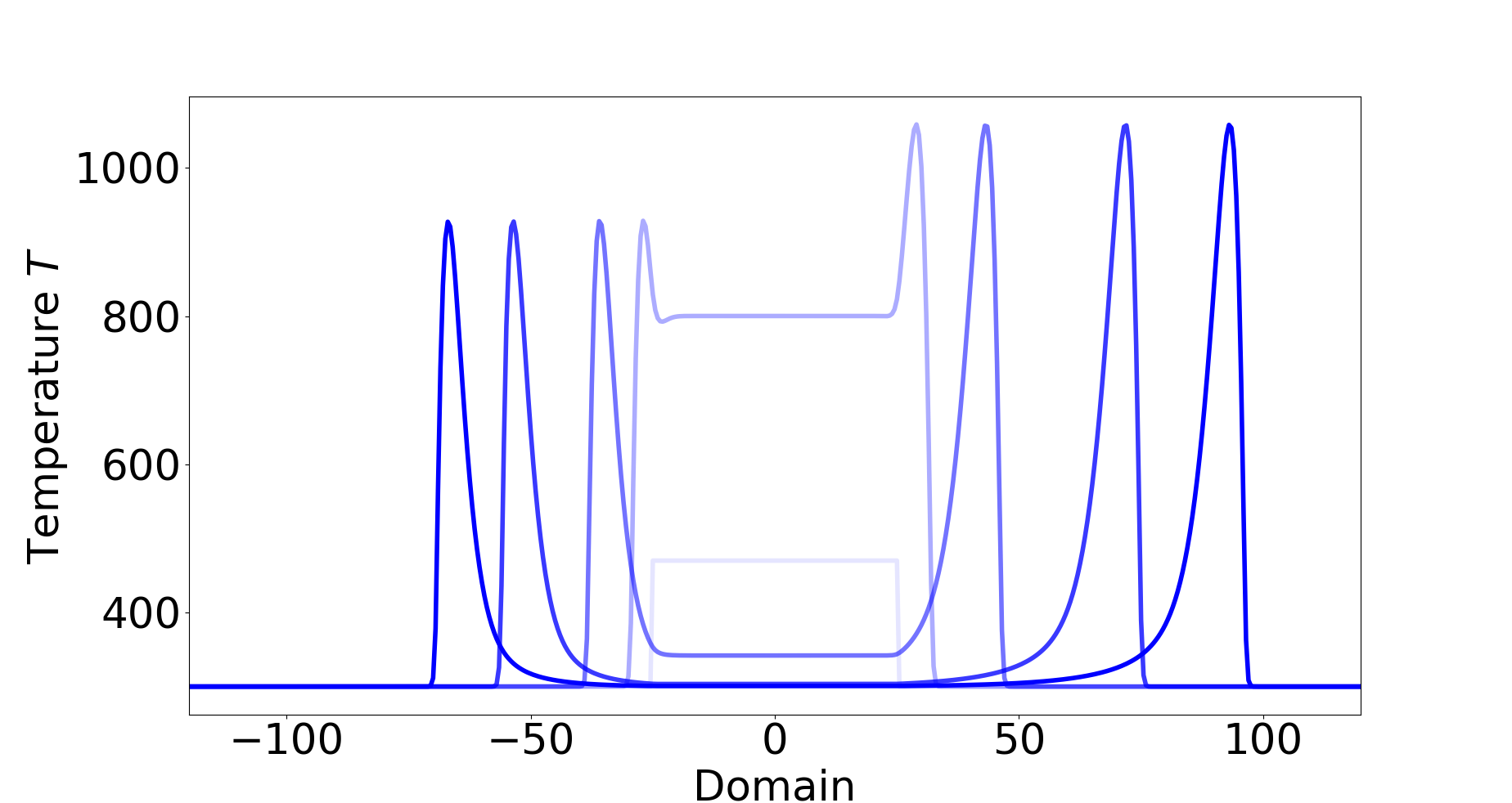}}
    \subfigure[Mild wind and TWs exist]{
    \includegraphics[width = 0.49\textwidth]{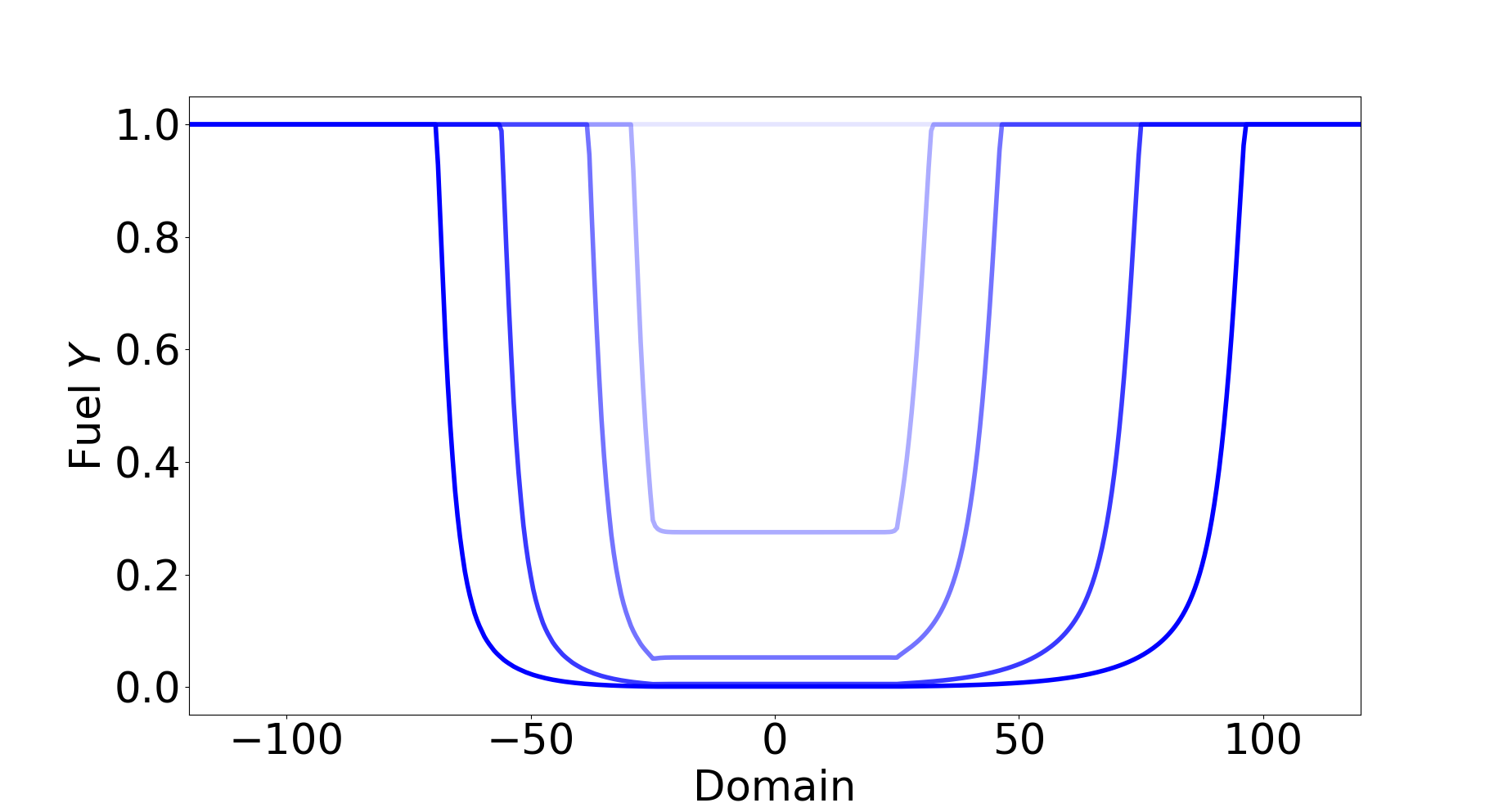}}
    \subfigure[Strong wind and TW exists]{
    \includegraphics[width = 0.49\textwidth,trim={5 0 5 0},clip]{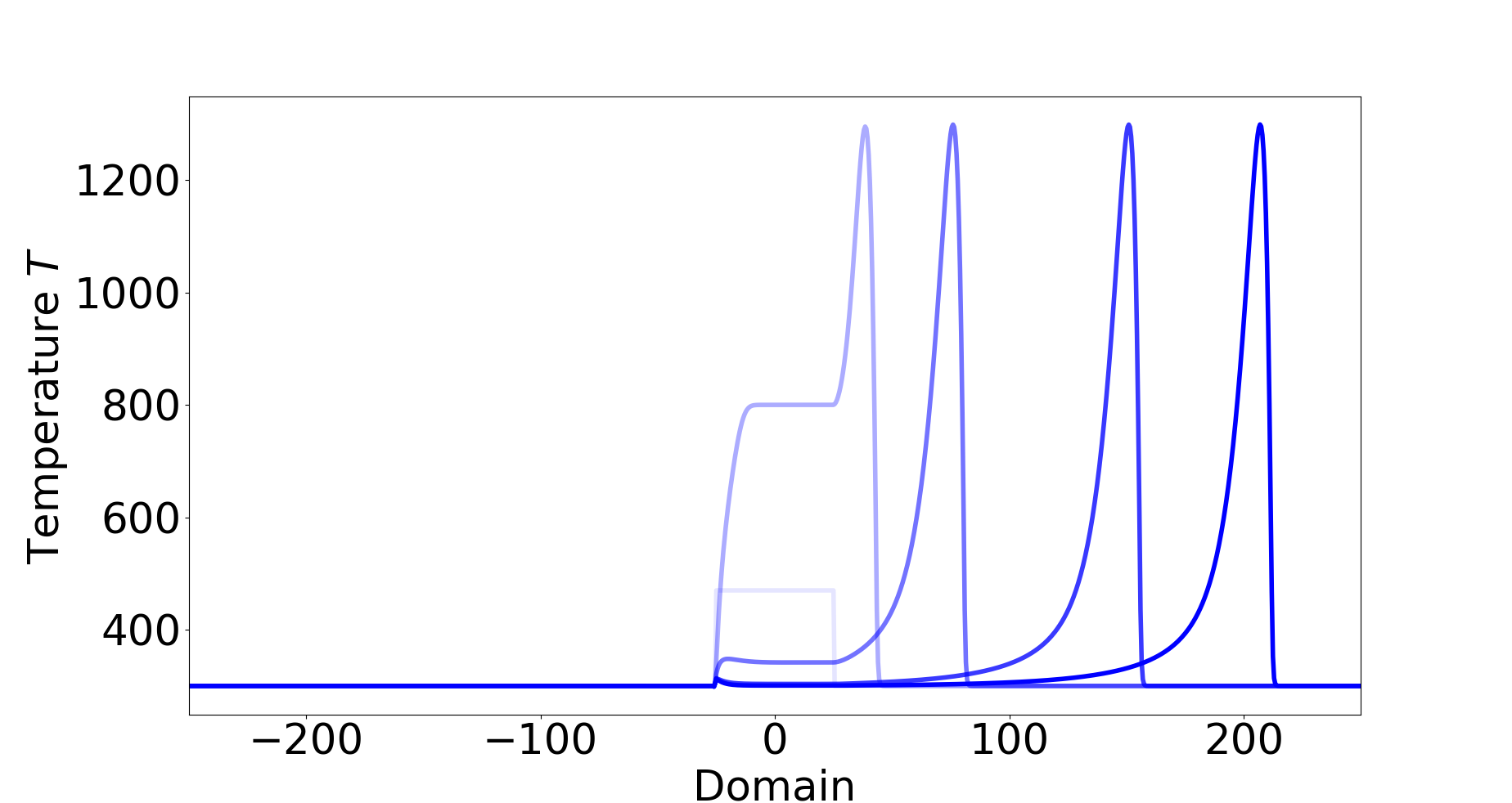}}
    \subfigure[Strong wind and TW exists]{
    \includegraphics[width = 0.49\textwidth]{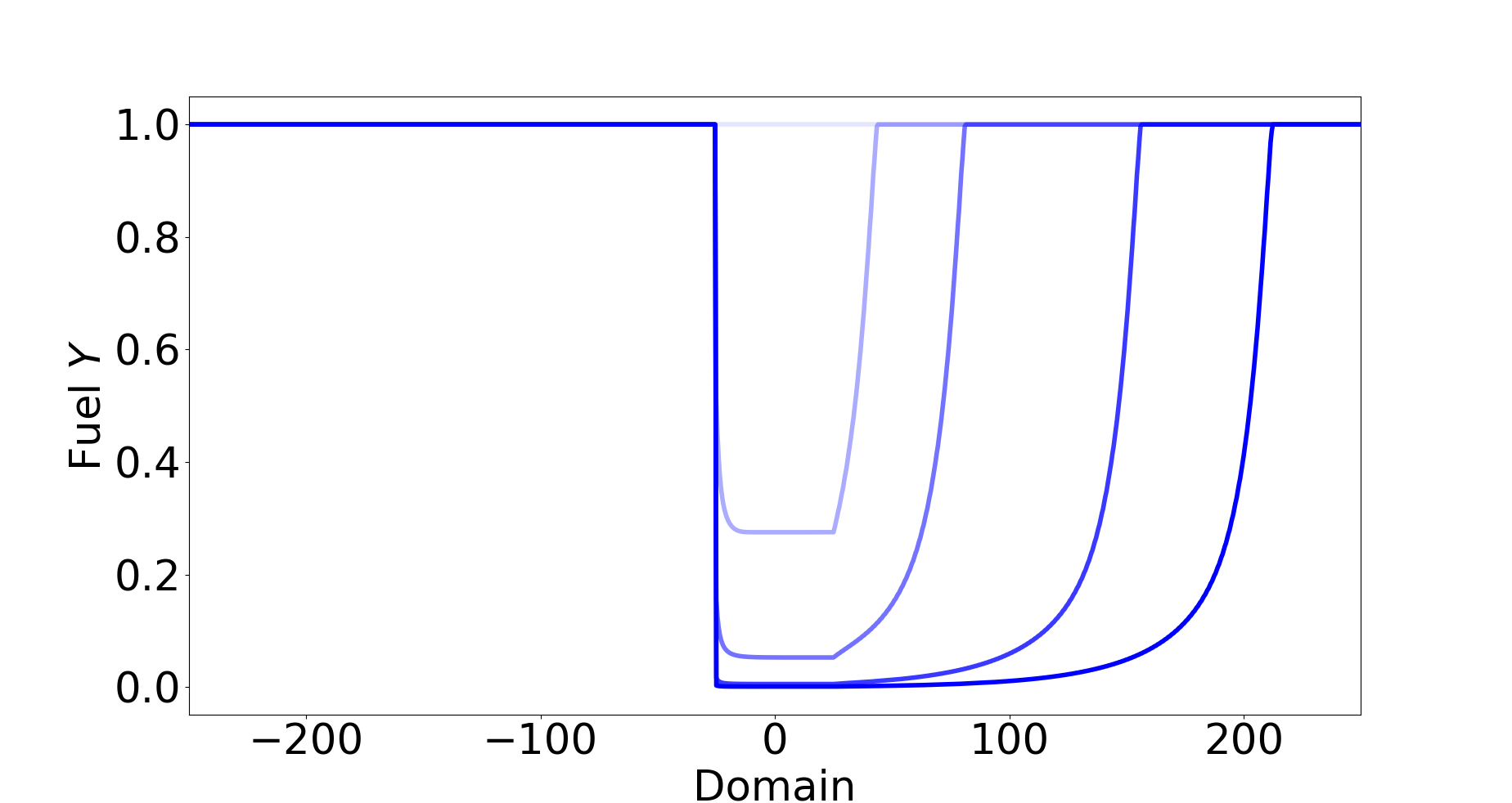}}
        \subfigure[No wind and no TW]{
    \includegraphics[width = 0.49\textwidth,trim={5 0 5 0},clip]{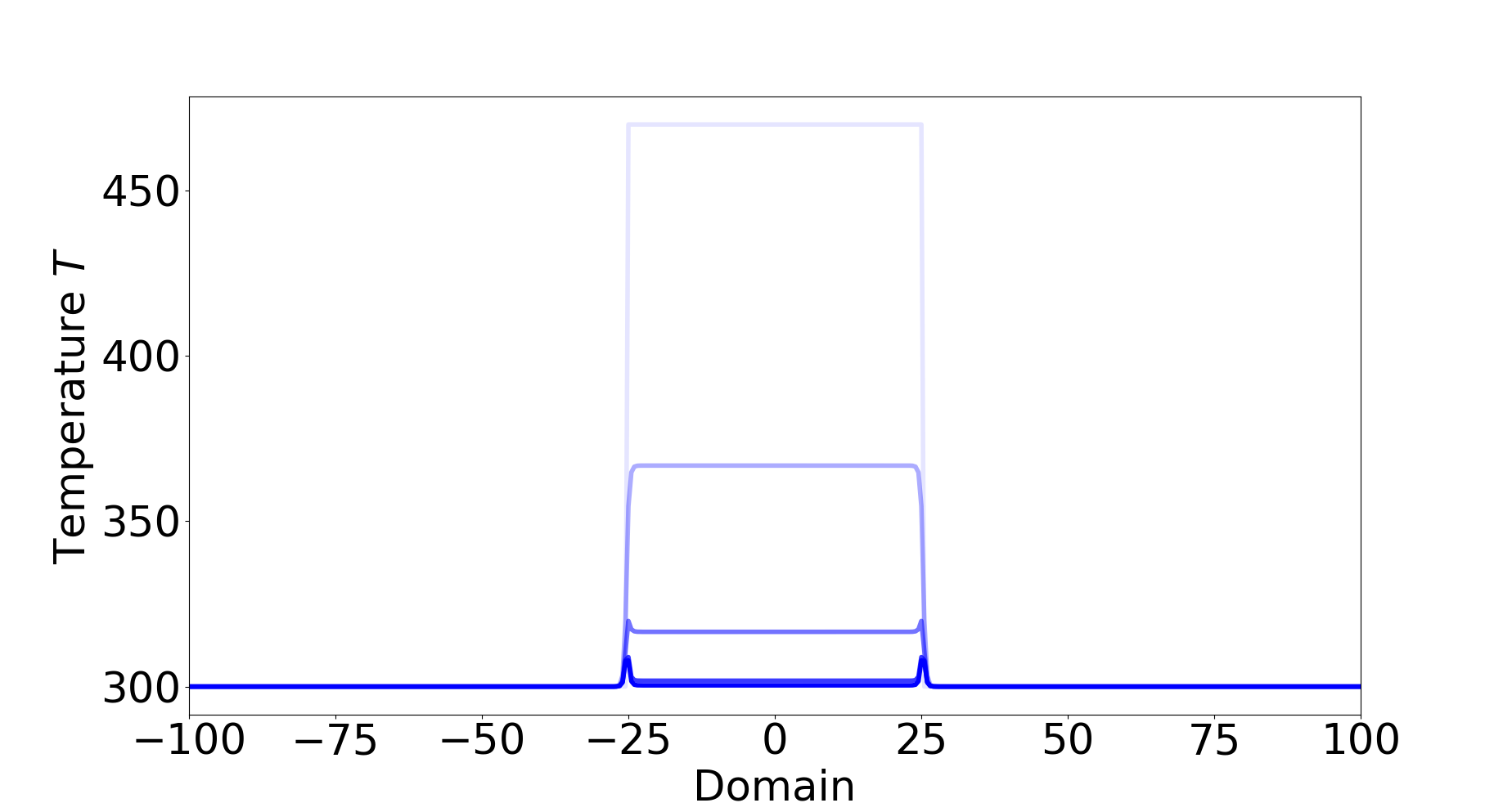}}
    \subfigure[No wind and no TW]{
    \includegraphics[width = 0.49\textwidth]{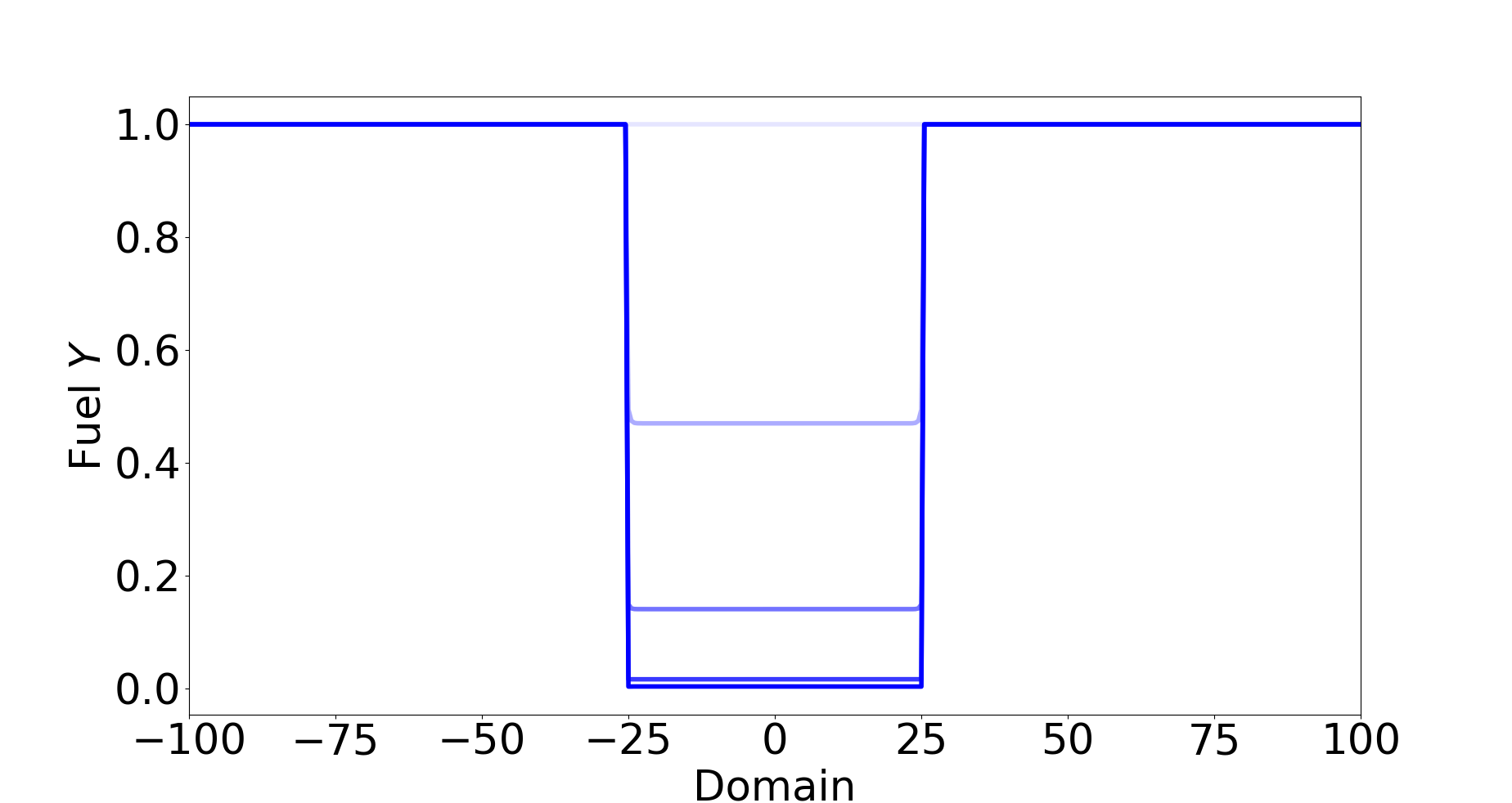}}
    \caption{The evolution of temperature $T$ (left column) and fuel $Y$ (right column) profiles are shown over (part of) the computational domain. The transparency of the profiles decreases as the time increases ($t \in \{0, 40, 120, 280, 400\}$). Here, (a)-(b) represents the no wind ($w=0 \,m/s$) case; (c)-(d) mild wind ($w=0.03m/s$); (e)-(f) strong wind ($w=0.3m/s$); and (g)-(h) no wind case where TWs do not exist (dimensional $h=19.85 kW/(m^3\cdot K)$, dimensionless $h=1.5$).}
    \label{fig:PDE_Eq3_exp}
    \vspace{-2em}
\end{figure}

Three main test cases are considered in this work for the TWs: (i) no wind $w=0m/s$, (ii) mild wind $w=0.03m/s$, and (iii) strong wind $w=0.3m/s$.  The evolution of temperature $T$ and the fuel concentration $Y$ profiles are shown in Figure \ref{fig:PDE_Eq3_exp} for these three test-cases. They show that the profiles immediately develop into TWs. For the no wind case, there are two symmetric TWs as expected. For the mild wind case, they are asymmetric with $c_++c_->0$; and for strong wind, there is only one TW with $c_+>0$ due to strong advection. The TWs speeds are calculated by linear regression on points $x(t)$ ($t>0$) for which $T(x(t),t)=\bar{T}$, see \Cref{tab:TW_speed}. The results support our claim in \Cref{sec:TW}. Furthermore, a fourth case is considered where no TWs exist even though there is no wind. It also agrees with our predictions, since the scaled convection constant $h>1$ in this case.

\begin{figure}[h]
    \centering
    \subfigure[No wind]{
    \includegraphics[width = 0.49\textwidth]{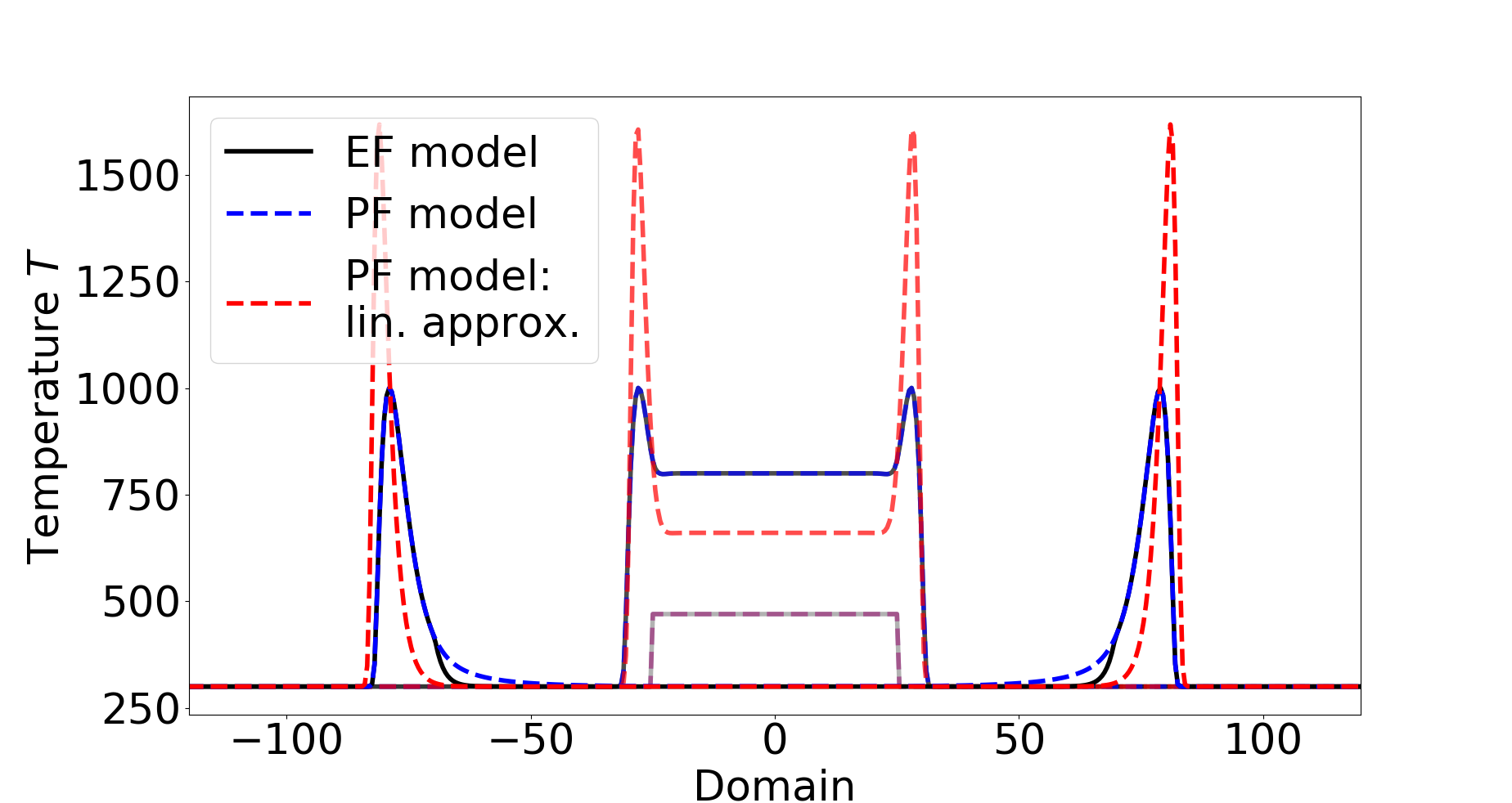}}
    \subfigure[No wind]{
    \includegraphics[width = 0.49\textwidth]{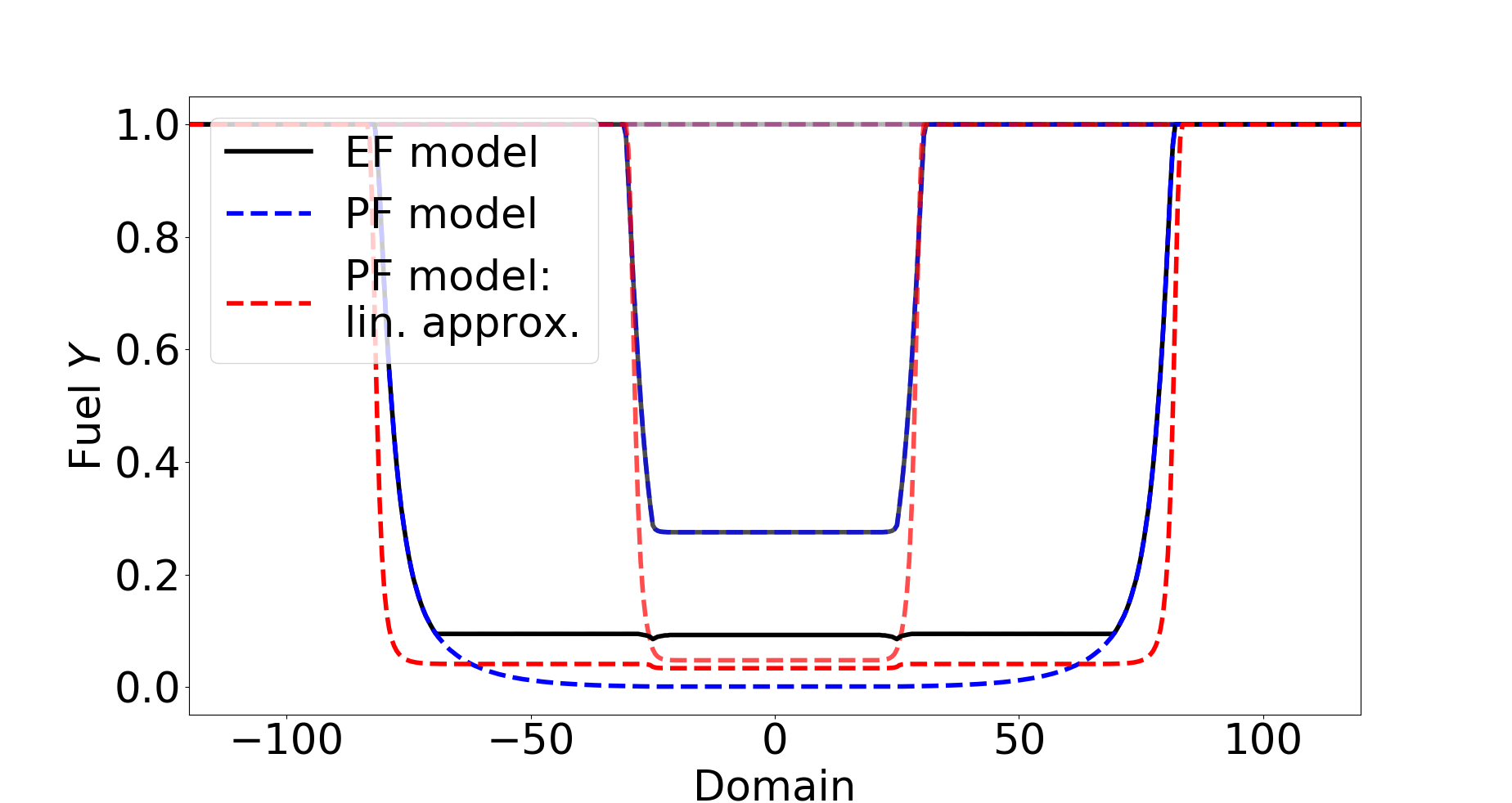}}
    \subfigure[Mild wind $(w = 0.03)$]{
    \includegraphics[width = 0.49\textwidth]{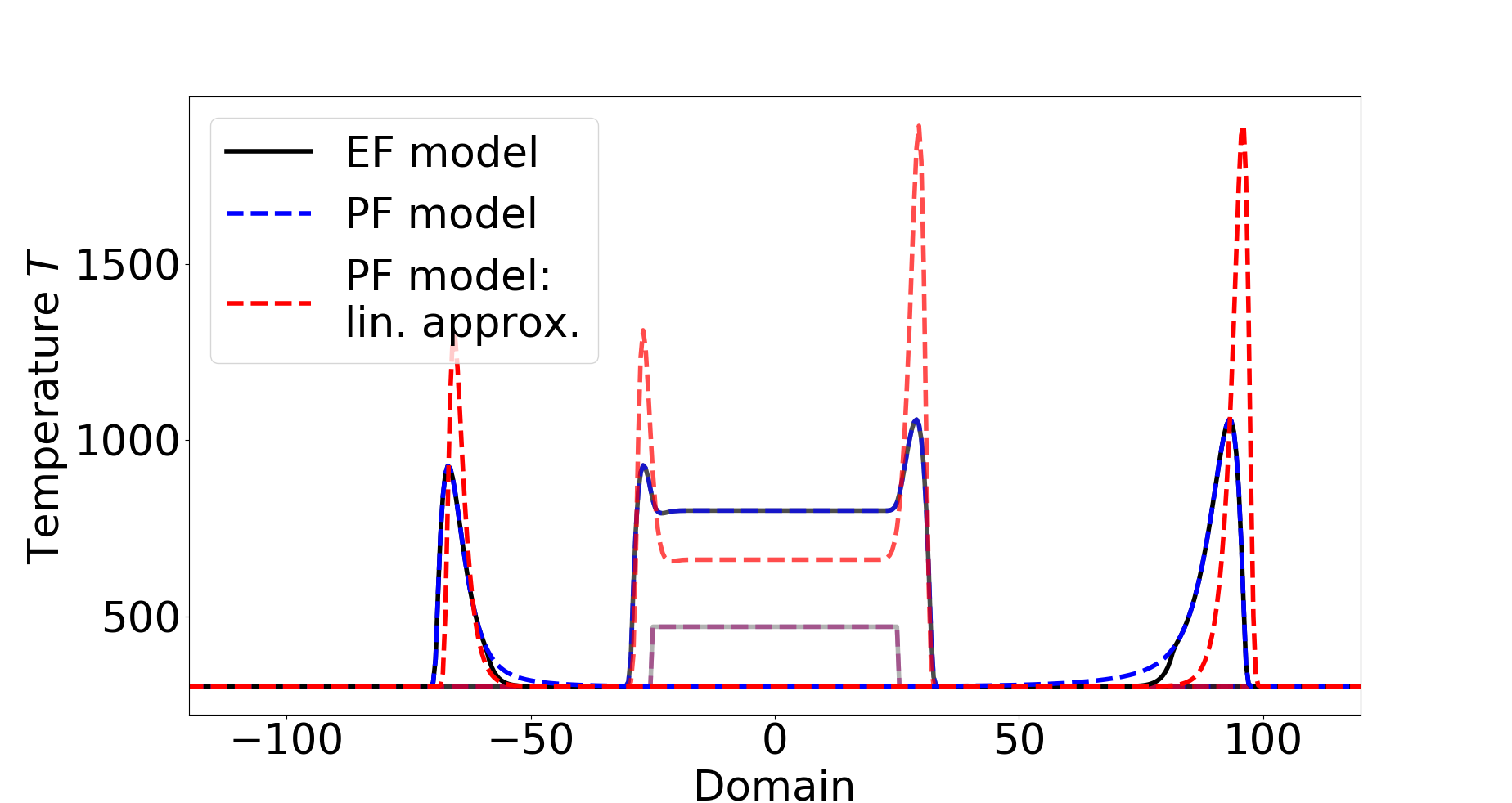}}
    \subfigure[Mild wind $(w = 0.03)$]{
    \includegraphics[width = 0.49\textwidth]{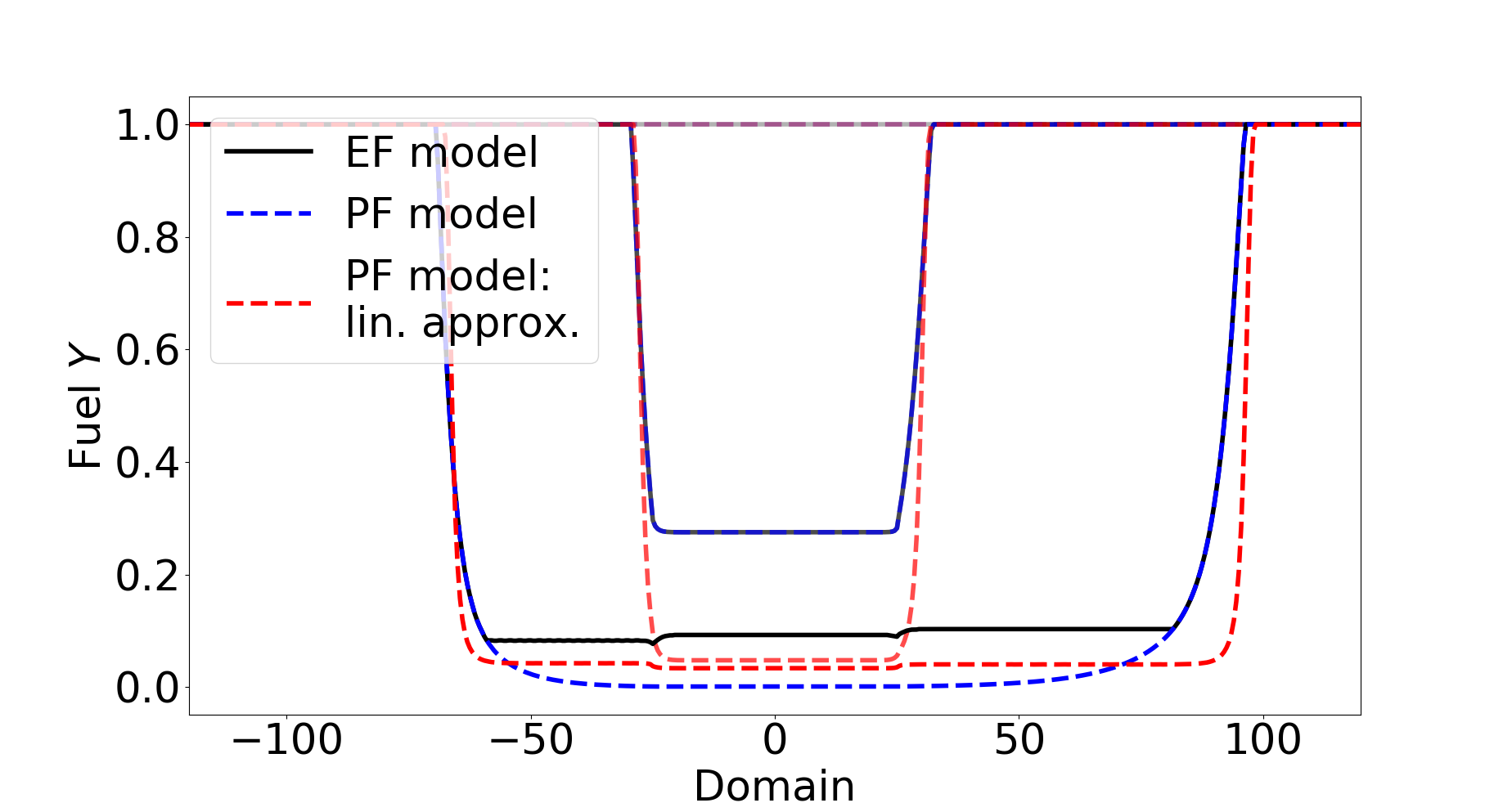}}
    \subfigure[Strong wind $(w = 0.3)$]{
    \includegraphics[width = 0.49\textwidth]{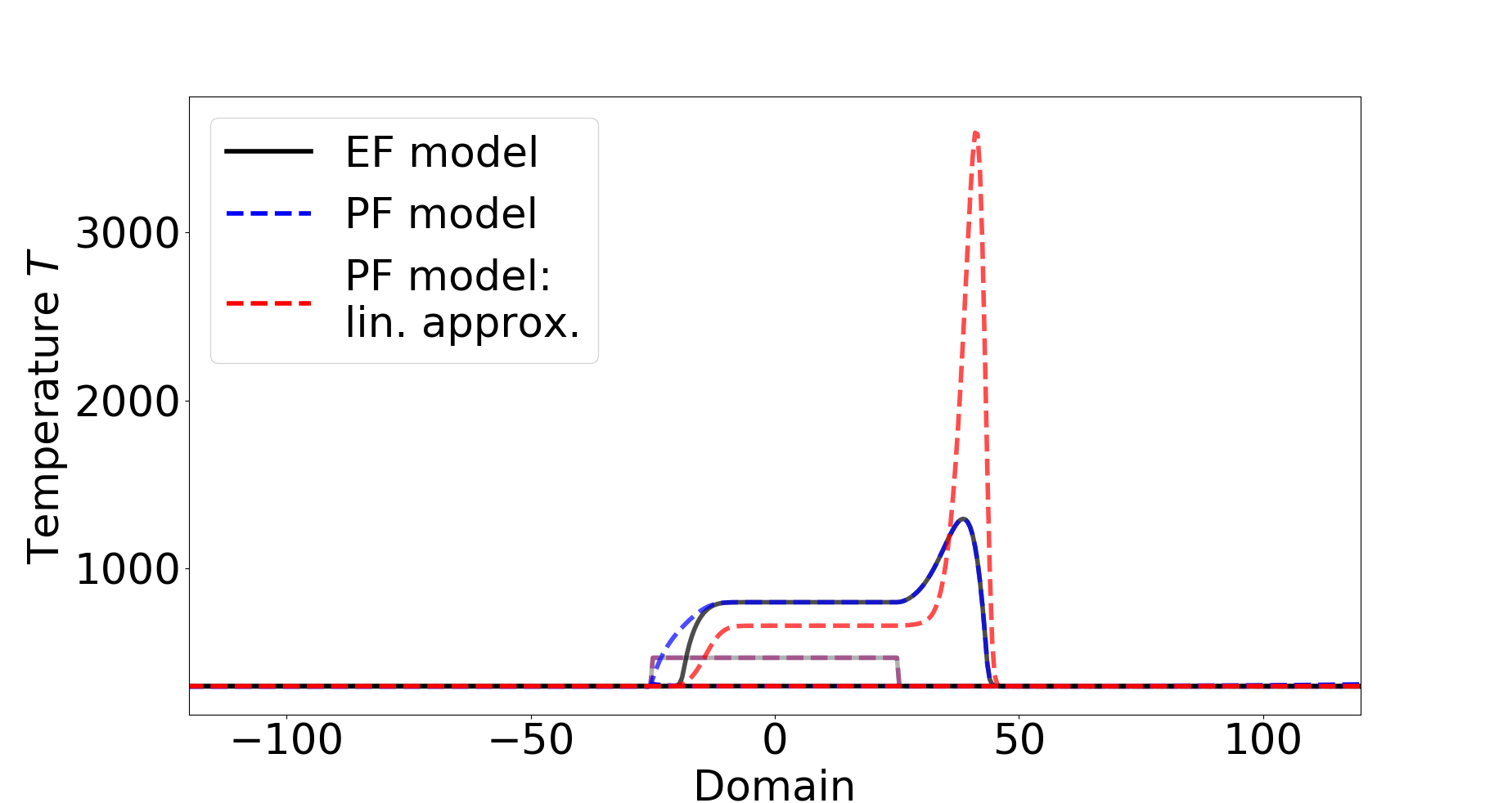}}
    \subfigure[Strong wind $(w = 0.3)$]{
    \includegraphics[width = 0.49\textwidth]{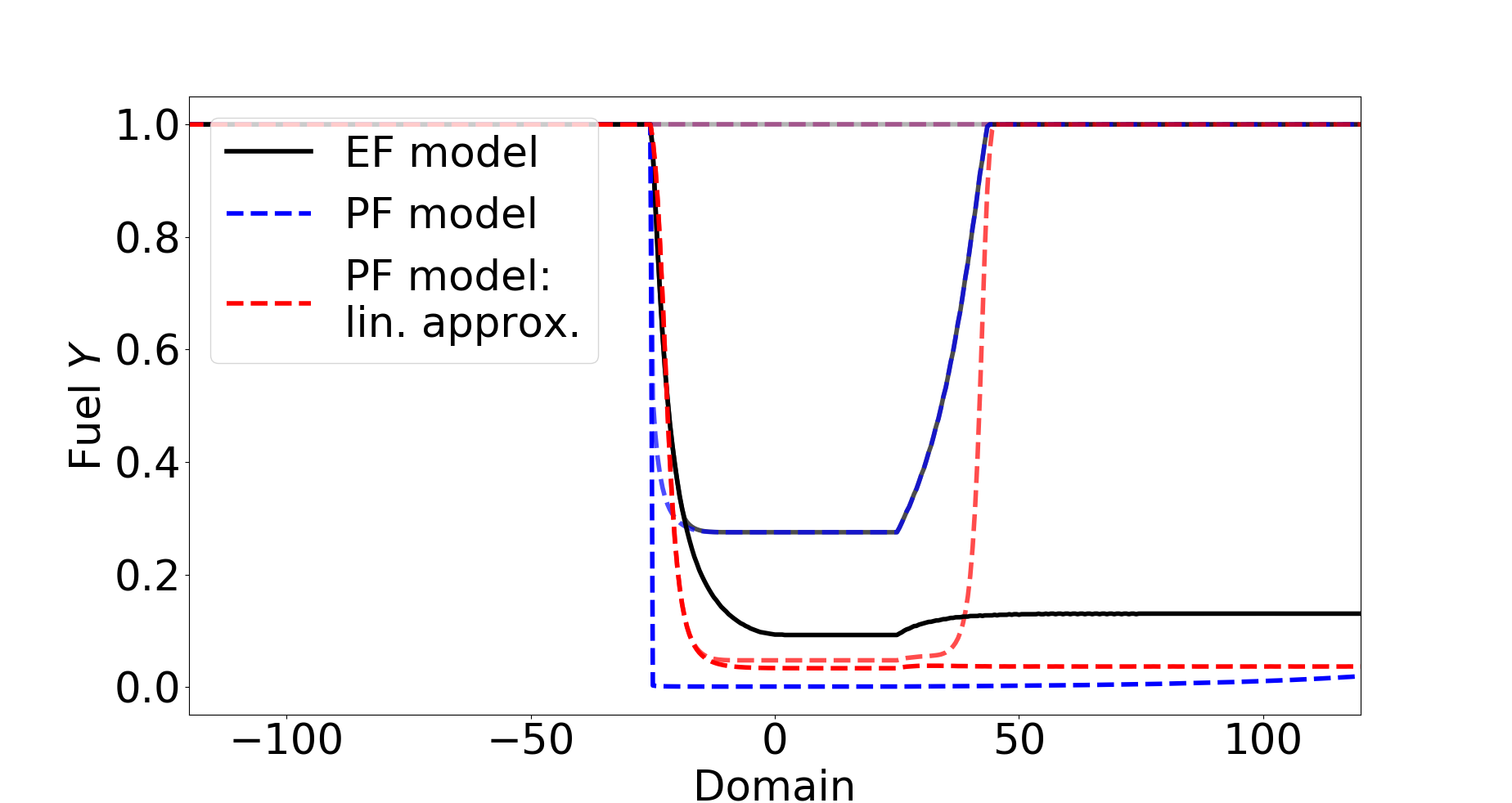}}
    \caption{The comparison between profiles of ADR-EF (black solid), ADR-PF (blue dashed), and the approximated (red dashed) models shown at various times ($t\in\{0,40,400\}$).}
    \label{fig:Model_comparison}
\end{figure}

Next, Figure \ref{fig:Model_comparison} compares three models mentioned in \Cref{sec:ADRmodel}, i.e. ADR-EF, ADR-PF, and ADR-PF with approximation \eqref{eq:linear_approx}, for the three test-cases. It is observed that the profiles for ADR-EF and ADR-PF models coincide until the tail region, where $T<\bar{T}$, is reached. Hence the existence of TWs and their speeds agree in both cases confirming the observations in \cite{weber1997} that the tail-region does not play a crucial role in this regard. However, the residual fuel $Y$ differs significantly in both cases which exemplifies why proper modeling is necessary for these cases. On the other hand, the approximate model results in much higher peak temperatures. This is because the exponential heat generation factor in \eqref{eq:new_psi2} plateaus at higher temperatures, whereas, the linear one in \eqref{eq:linear_approx} does not. However, this approximate model still manages to predict the existence and speed of the TWs in all three cases revealing its usefulness as an analytical tool, see \Cref{tab:TW_speed}.


\begin{figure}[h]
\subfigure[no wind $w=0$]{    \includegraphics[width=0.32\textwidth]{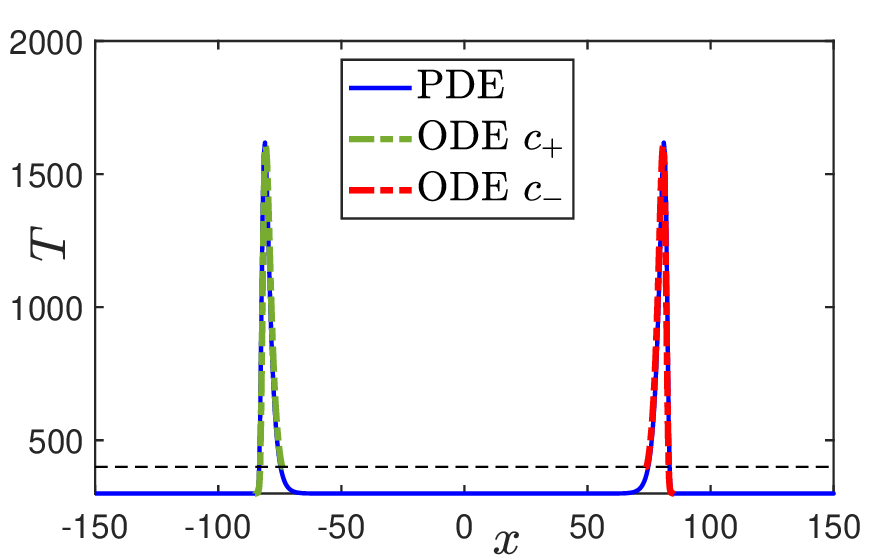}}
\subfigure[mild wind $w=0.03m/s$]{    \includegraphics[width=0.32\textwidth]{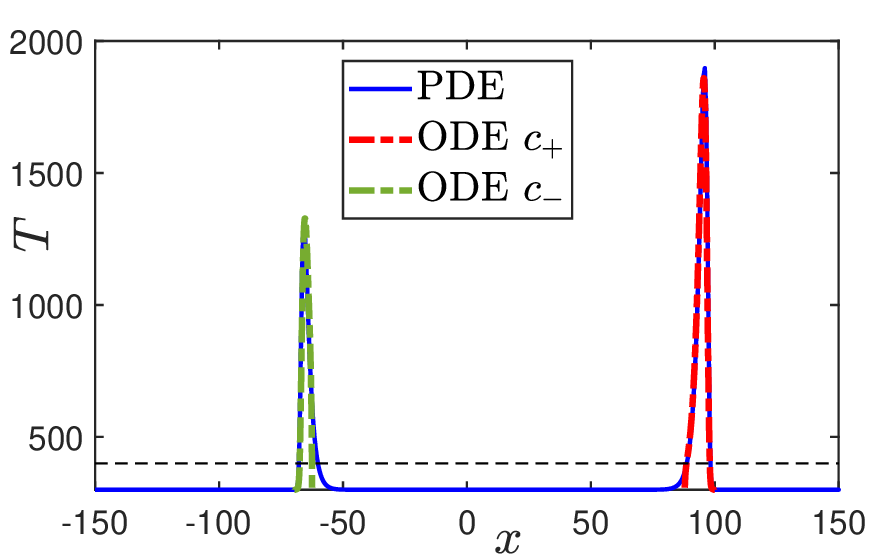}}
\subfigure[strong wind $w=0.3m/s$]{    \includegraphics[width=0.32\textwidth]{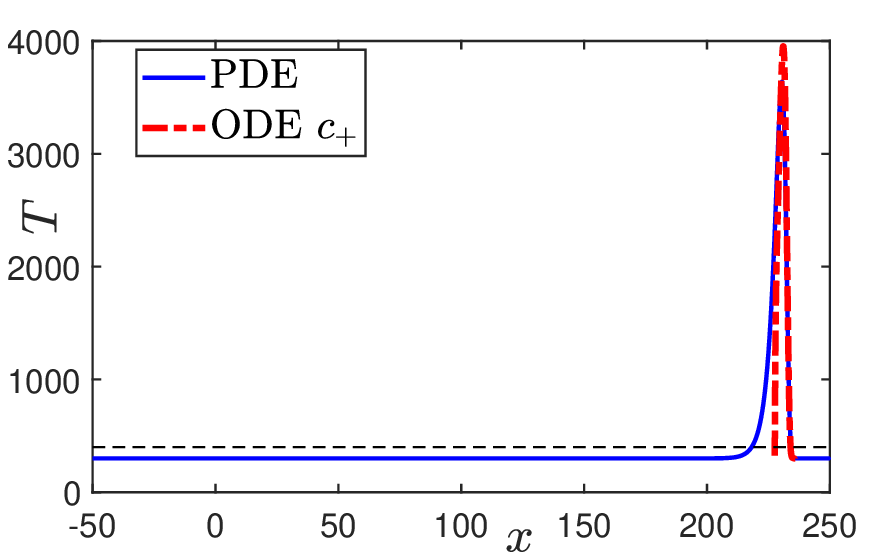}}
    \caption{Comparison of temperature profiles between the PDE and the ODE simulations of the approximated ADR-PF model.}
    \label{fig:pde-ode}
    \vspace{-1.5em}
\end{figure}

Next, we compare the PDE results of the approximated ADR-PF model \eqref{eq:original_model}, \eqref{eq:new_psi2}--\eqref{eq:linear_approx} with solutions of the dynamical system \eqref{eq:uvz_ODE_xi} with conditions \eqref{eq:BC_ode}. To find these ODE solutions, we use 
\vspace{-.5em}
\begin{algorithm}[A shooting ODE solver]\label{algo:shooting}
For a given $c>0$, let $(u_c,v_c,z_c)$ denote the solution of \eqref{eq:uvz_ODE_xi} connecting to $(0,1,0)$ as $\xi\to \infty$ with $u_c(0)=1$. Define,
$d_c:= \inf\{u_c(\xi)^2 + (hz_c(\xi)-\gamma^{-1}(1-v_c(\xi)))^2: \xi<0 \}$
as the closest approach of $(u_c,v_c,z_c)$ to another equilibrium point of System \eqref{eq:uvz_ODE_xi}. Then assuming that the unique minimizer of $d_c$ is in an interval $[\Check{c},\hat{c}]$ for some $0<\Check{c}<\hat{c}$, we employ a line-search algorithm to find $c=c_+$ such that $d_c\approx 0$.
\vspace{-.5em}
\end{algorithm}
For finding $c=c_-<0$, the same process is repeated taking $w'=-w$ due to symmetry.
For the actual ODE simulations, a fourth-order Runge-Kutta scheme has been used to compute backward from $(10^{-3},1,0)$ with $\Delta \xi=10^{-3}$, and $(\Check{c},\hat{c})=(0.1,1.0)$. \Cref{fig:pde-ode} shows the comparison of temperature profiles resulting from the PDE (approximated ADR-PF) and the ODE solutions for the three test cases. The profiles are a close match to each other until the orbits of the ODE solutions are near the second equilibrium point. \Cref{tab:TW_speed} shows a comparison between TW speeds for the different solvers which differ from each other by at most $4\%$. Overall, we find the matching to be excellent.

\begin{table}[h!]
    \centering
    \begin{tabular}{|c|c|c|c|c|}
        \hline
         \textbf{Case} & \textbf{Wind speed $w$} & \textbf{ADR-PF $(c_-, c_+)$} & \textbf{Lin. Approx. $(c_-, c_+)$} & \textbf{ODE $(c_-, c_+)$}   \\
         \hline
         No wind & $0$ & $(-0.1421, 0.1421)$ & $(-0.1585, 0.1585)$ & $(-0.15196,0.15196)$\\
         \hline
         Mild wind & $0.03$ & $(-0.1112,0.1781)$ & $(-0.1079, 0.1843)$ &$(-0.116358,0.18585)$   \\
         \hline
        Strong wind & $0.3$ & $({\rm NA}, 0.4686)$ & $({\rm NA}, 0.5264)$ &$({\rm NA},0.4702196)$ \\
        \hline
    \end{tabular}
    \caption{Comparison of wave-speeds of the numerical solutions of the ADR-PF model (Eqs. \eqref{eq:original_model}, \eqref{eq:new_psi2}, the approximated ADR-PF model (Eqs. \eqref{eq:original_model}, \eqref{eq:new_psi2}--\eqref{eq:linear_approx}), and  \Cref{algo:shooting} (ODE solver). All three test cases (no, mild, and strong wind) are considered.}
    \label{tab:TW_speed}
    \vspace{-2em}
\end{table}

\section{Discussion and future directions}
\vspace{-1em}
We studied an existing advection--diffusion--reaction model for wildfire propagation and proposed an alteration in which fire persists behind the front. A one-dimensional version of the model was considered with an approximation used for estimating the Arrhenius factor. Using travelling wave ansatz, we reduced the system to a dynamical system with boundary values prescribed. This system is easier to analyze and solve rapidly using a shooting algorithm (\Cref{algo:shooting}). It was also hypothesized that the travelling waves emanate in both directions under mild wind conditions, but only in the direction of wind provided the wind speed is high enough. Our predictions were corroborated by numerical results. 

The algorithm can be used to quickly predict the wildfire propagation speed under different parametric conditions to prevent disasters. The actual proof of the existence of the travelling waves will be elaborated in a future work.
\vspace{-1em}

\begin{acknowledgement}
CR acknowledges funding 
by the Seed Funding Programme of Technische Universit\"at Braunschweig, 2022 Interdisciplinary Collaboration: Strengthening Interdisciplinarity
-- Expanding Research Collaboration. KM was supported by Fonds Wetenschappelijk Onderzoek (FWO) through the Junior Postdoctoral Fellowship.
\end{acknowledgement}
%



\vspace{-2em}
%
%
\bibliographystyle{unsrtnat}

\begin{thebibliography}{99.}%
%
%
\bibitem{asensio_2002} Asensio, M. I.,  Ferragut, L.: On a wildland fire model with radiation.  International Journal
for Numerical Methods in Engineering \textbf{54}, 137--157 (2002). doi: 10.1002/nme.420.

\bibitem{reisch_2023} Reisch, C., Navas-Montilla, A., {\"O}zgen-Xian, I.: Analytical and numerical insights into wildfire dynamics: Exploring the advection-diffusion-reaction model. arXiv:2307.16174 (2023)

\bibitem{senande_2022} Senande-Rivera, M., Insua-Costa, D., Miguez-Macho, G.: Spatial and temporal expansion of global wildland fire activity in response to climate change. Nat Commun \textbf{13}, 1208 (2022). doi:10.1038/s41467-022-28835-2

\bibitem{burger_2020b} B{\"u}rger, R., Gavi\'{a}n, E.,  Inzunza, D., Mulet, P., Villada, L. M.: Implicit-Explicit Methods for a Convection-Diffusion-Reaction Model of the Propagation of Forest Fires. Mathematics \textbf{8}, 1034 (2020) doi: 10.3390/math8061034

\bibitem{weber1991} Weber, R.O.: Toward a Comprehensive Wildfire Spread Model.  Int. J. Wildland Fire \textbf{1}, 245-248 (1991). doi: 10.1071/WF9910245

\bibitem{VolpertTW1994}
Volpert, A.I., Volpert, V.A., Volpert, V.A.: Traveling Wave Solutions of Parabolic Systems. Translations of Mathematical Monographs 140, AMS (1994). doi: 10.1090/mmono/140

\bibitem{weber1997} Weber R.O., Mercer G.N., Sidhu H.S., Gray B.F.: Combustion waves for gases ($Le = 1$) and solids ($Le \to \infty$). Proc. R. Soc. Lond. A. \textbf{453}, 1105–1118 (1997). doi: 10.1098/rspa.1997.0062


\bibitem{mandel2008} Mandel, J., Bennethum, L.S., Beezley, J.D., Coen, J.L.,  Douglas, C.C., Kim, M., Vodacek, A.: A wildland fire model with data assimilation. Mathematics and Computers in Simulation \textbf{79}(3), 584-606 (2008). doi: 10.1016/j.matcom.2008.03.015.


\bibitem{burger_2020} B{\"u}rger, R., Gavi\'{a}n, E.,  Inzunza, D., Mulet, P., Villada, L. M.:Exploring a convection-diffusion-reaction model of the propagation of forest fires: Computation of risk maps for heterogeneous environments. Mathematics \textbf{8}, 1674 (2020) doi: 10.3390/math8101674


\bibitem{babak_2009} Babak, P., Bourlioux, A., Hillen, T.: The Effect of Wind on the Propagation of an Idealized Forest
Fire. SIAM Journal on Applied Mathematics \textbf{70}(4) 1364--1388 (2009) doi: 10.1137/080727166.



\bibitem{MitraTW2023}
Mitra, K., Hughes, J.M., Sonner, S., Eberl, H.J., Dockery, J.D.: Travelling Waves in a PDE–ODE Coupled Model of Cellulolytic Biofilms with Nonlinear Diffusion. J Dyn Diff Equat (2023). https://doi.org/10.1007/s10884-022-10240-4




\end{thebibliography}
%

\end{document}